\documentclass[12pt]{amsart}
\usepackage{amsfonts,amsmath,amsthm,amscd,amssymb,latexsym,cite,verbatim,texdraw,floatflt,
caption2}
\RequirePackage{hyperref}

\usepackage[a4paper]{geometry}

\begin{document}

\title{  Coarse structures  on groups defined by conjugations }

\author{ Igor Protasov, Ksenia Protasova}

\maketitle
\vskip 5pt

{\bf Abstract.} For a group $G$, we denote by 
$\stackrel{\leftrightarrow}{G}$ the coarse space on $G$ endowed with the coarse structure with the base $\{\{ (x,y)\in G\times G: y\in x^F \} : F \in [G]^{<\omega} \}$, 
$x^F = \{z^{-1} xz : z\in F \}$. Our goal is to explore interplays between algebraic properties of $G$ and asymptotic properties of 
$\stackrel{\leftrightarrow}{G}$. 
In particular, we show that $asdim \ \stackrel{\leftrightarrow}{G} = 0$ if and only if  $G / Z_G$ is locally finite, $Z_G$
is the center of $G$. For an infinite group $G$, the coarse space of subgroups of $G$ is discrete if and only if $G$  is a Dedekind group. 

\vspace{6 mm}

 20E45, 54D80

\vspace{3 mm}

Keywords:  coarse structure defined by conjugations, cellularity, FC-group, ultrafilter.


\section{Introduction }

Given a set $X$, a family $\mathcal{E}$ of subsets of $X\times X$ is called a {\it coarse structure }  on $X$  if
\vskip 7pt

\begin{itemize}
\item{}   each $E\in \mathcal{E}$  contains the diagonal  $\bigtriangleup _{X}$,
$\bigtriangleup _{X}= \{(x,x)\in X: x\in X\}$;
\vskip 5pt

\item{}  if  $E$, $E^{\prime} \in \mathcal{E}$ then $E\circ E^{\prime}\in\mathcal{E}$ and
$E^{-1}\in \mathcal{E}$,   where    $E\circ E^{\prime}=\{(x,y): \exists z((x,z) \in  E,  \   \ (z, y)\in E^{\prime})\}$,   $E^{-1}=\{(y,x): (x,y)\in E\}$;
\vskip 5pt

\item{} if $E\in\mathcal{E}$ and $\bigtriangleup_{X}\subseteq E^{\prime}\subseteq E  $   then
$E^{\prime}\in \mathcal{E}$;
\vskip 5pt

\end{itemize}
\vskip 7pt

A subfamily $\mathcal{E}^{\prime} \subseteq \mathcal{E}$  is called a
{\it base} for $\mathcal{E}$  if,
 for every $E\in \mathcal{E}$, there exists
  $E^{\prime}\in \mathcal{E}^{\prime}$  such  that
  $E\subseteq E ^{\prime}$.
For $x\in X$,  $A\subseteq  X$  and
$E\in \mathcal{E}$, we denote
$$E[x]= \{y\in X: (x,y) \in E\},
 \   E [A] = \bigcup_{a\in A}   \    
  E[a], \   \     
   E_A   [x]  = E[x]\cap A$$
 and say that  $E[x]$
  and $E[A]$
   are {\it balls of radius $E$
   around} $x$  and $A$.

\vskip 10pt
The pair $(X,\mathcal{E})$ is called a {\it coarse space}  
 \cite{b13} 
 or a {\it ballean}
\cite{b10}, 
\cite{b12}. 


 A coarse space $(X,\mathcal{E})$
 is  called {\it  finitary, } if for
 each 
 $E\in   \mathcal{E} $, there exists a natural number $n$ such that $|E[x]|< n$  for each $x\in X$.

Let $G $ be a group of permutations of a set $X$. 
We denote by $X_G$ the set $X$ endowed with the coarse structure with the base  
$$\{\{ (x,gx): g\in F \}:  F\in [G]^{<\omega}, \ id\in F\} .$$

By  [7,  Theorem 1], for every finitary coarse structure 
 $(X, \mathcal{E})$,  there exists a group  $G$ of permutations of $X$ such that   $(X, \mathcal{E}) = X_G. $ 
For more general results and applications  see 
\cite{b8} and the survey \cite{b9}.

\vskip 7pt

Let  $(X, \mathcal{E})$ be a  coarse space. We define an equivalence $ \sim $ on $X$ by $x \sim y$ if and only if there exists $E\in \mathcal{E}$ such that $y\in E[x]$, so $X$ is a disjoint union of connected components. If there is only one connected component then $(X, \mathcal{E})$ is called connected. 

Now let $G$ be a group. For $x,g\in G$  and $F\subseteq G$, we denote $x^g = g^{-1} x g$, $x^F = \{ x^{y}: y\in F\}$, 
$F^g  = \{ y^{g}: y\in F\}$.

We denote by  
$\stackrel{\leftrightarrow}{G}$ the coarse structure on $G$ endowed with the coarse structure with the base 
$\{\{ (x,y)\in G\times G: y\in x^F \}:  F\in [G]^{<\omega} \} .$
Evidently, each connected component $A$ of 
$\stackrel{\leftrightarrow}{G}$ is of the form $a^G$, $a\in A$.

We endow $G$ with the discrete topology and identify the 
Stone-$\check{C}$ech compactification $\beta G$ of $G$ with the set of all ultrafilters on $G$. 
For $A\subseteq G$, $\bar{A}$ denotes the set $\{ p\in \beta G: A\in p \}$ and the family 
$\{ \bar{A}: A\subseteq G \}$ forms a base for open sets of $\beta G$. The family of all free ultrafilters on $G$ is denoted by $G^\ast$. 
By the universal property of $\beta G$, every mapping $f: G\rightarrow K$, $K$ is a compact Hausdorff space, 
can be extended  to
the continuous mapping $f^\beta : \beta G\rightarrow K $.

The action $G$ on $G$ by conjugations extends to the action $G$  on $\beta G: $ if $g\in G$,  $p\in \beta G$ then $p^g = \{g^{-1} Pg: P\in g\} $.
We use this dynamical approach to the conjugacy in groups initiated in \cite{b11}.

\vskip 7pt

In section 2 and 3, we characterize groups $G$ such that the coarse space 
$\stackrel{\leftrightarrow}{G}$ is discrete, $n$-discrete and cellular. In section 4, we show that every finitary coarse space admits an asymorphic embedding to $\stackrel{\leftrightarrow}{G}$  for an appropriate choice of a group $G$. In section 5, we characterize groups with discrete space of subgroups.
We conclude with section 6 on the direct union of connected components of $\stackrel{\leftrightarrow}{G}$.

\section{Discreteness}

Let $(X, \mathcal{E})$  be a coarse space. We say that a subset $B$ of $X$ is {\it bounded} if there exist a finite 
subset $F$ of $X$
and $E\in \mathcal{E}$ such that $B\subseteq E[F]$ and note that the family of all bounded subset of $X$  is a bornology, i.e. an ideal in the Boolean algebra of subsets of $X$ containing all finite subsets. 

We say that a subset $A $ of $X$ is  

\vskip 7pt

\begin{itemize}

\item{}  {\it discrete}  if, for every   $E\in \mathcal{E}$, there exists a bounded subset $B$ of $X$ such that $E_A [a]= \{a\}$ for each $a\in A\setminus B$;

\vskip 7pt

\item{}  n-{\it discrete}, $n\in \mathbb{N}$  if, for every $E\in \mathcal{E}$, there exists a bounded subset $B$  of $X$ such that $|E_A [a]|\leq n$ for each $a\in A\setminus B$.

\end{itemize}

\vskip 7pt

{\bf Theorem 1.} {\it For an infinite group $G$, the following conditions are equivalent 

\vskip 10pt

(i)  $G$ is Abelian;

\vskip 10pt

(ii)  $p^G = \{p\}$  for each $p\in G^\ast$;

\vskip 10pt

(iii)  $\stackrel{\leftrightarrow}{G} $   is discrete.

\vskip 10pt

Proof.}
The equivalence $(i) \Longleftrightarrow  \ (ii)$ is proved in 
[11, Proposition 1.1], (i) $\Longrightarrow$  (iii) is evident.
 \vskip 10pt

$(iii) \Longrightarrow  \ (ii)$. We assume that $p^x \neq p$
for  some $p\in G^\ast$, $x\in G$ and pick $P\in p$ such that $P^x \cap P = \emptyset$. Let $B$ be a finite subset of $X$. We take $a\in P\setminus B$  and note that $a^x \neq a$ so $\stackrel{\leftrightarrow}{G} $ is not discrete.  $ \ \Box$

\vskip 10pt

{\bf Theorem 2.} {\it For a group $G$, the following conditions are equivalent 

\vskip 10pt

(i)  $p^G$ is finite for each $p\in G^\ast$;

\vskip 10pt

(ii) there exists a natural number $n$ such that 
 $|p^G| \leq  n$  for each $p\in G^\ast$;

\vskip 10pt

(iii)  there exists a natural number $m$ such that 
$|a^G|\leq m$ for each $a\in G^\ast$;

\vskip 10pt

(iv) the commutant $[G, G]$ of $G$ is finite.

\vskip 10pt

Proof.} See Theorem 3.1 in \cite{b11}. $ \ \Box$

\vskip 10pt

{\bf Theorem 3.} {\it Given a group $G$, the 
coarse space $\stackrel{\leftrightarrow}{G} $ is 
$n$-discrete for some $n\in \mathbb{N}$ if and only if $[G, G]$ is finite.

\vskip 10pt

Proof.} We assume that $\stackrel{\leftrightarrow}{G} $ is 
$n$-discrete and show that $[G, G]$ is finite.
To apply Theorem 2, it suffices to prove that
 $|p^G|\leq n$ for each $p\in G^\ast$.

We assume the contrary: 
there exists 
$p\in G^\ast $ and $g_1, \dots , g_{n+1} \in G$ such that the ultrafilters $p^{g_1}, \dots , p^ {g_{n+1}}$ are distinct. 
We choose $P\in p$ such that the subsets $P^{g_1}, \dots , P^ {g_{n+1}}$ 
are pairwise disjoint. Given an arbitrary bounded subset $B$ of $G$, we pick $a\in P\setminus B$.
Then $a^{g_1}, \dots , a^ {g_{n+1}}$ are distinct so 
$\stackrel{\leftrightarrow}{G}$
is not $n$-discrete. 

On the other hand, if $[G, G]$ is finite then there exists $m\in \mathbb{N}$ such that $|a^G|\leq m$ for each $a\in G$, see  Theorem 2$(iii)$. 
$ \  \Box$

\vskip 10pt

We recall that $G$ is an {\it $FC$-group} if the set $a^G$ is finite for each $a\in G$. Clearly, $G$ is an $FC$-group if and only if each connected component of 
$\stackrel{\leftrightarrow}{G} $ 
is bounded.

\vskip 10pt


\vskip 10pt

We note that each connected component of 
$\stackrel{\leftrightarrow}{G} $ is discrete if and only if every element $g\in G$ centralizes all but finitely many elements of each conjugacy  class. 

In the initial version of this paper, we asked whether $G$ is an $FC$-group provided that each connected component of $\stackrel{\leftrightarrow}{G} $  is discrete?
G. Bergman answered this question negatively. 

\vskip 10pt

{\bf Theorem 4.} {\it There exists a group $G$ such that every element of $G$ centralizes all but finitely many element of each conjugacy class and 
$g^G$ is infinite  for each nonindentily element $g\in G$. 
\vskip 10pt

Proof.} We follow the original Bergman's exposition. 
\vskip 10pt

{\it Claim 1.} Suppose $X$ is a metric space such that, for every $x\in X$ and constant $C>0$, the number of elements of $X$ within distance $\leq C$ of $x$ is finite. 
Suppose also that $X$ has a group
$G$
 of distance-preserving permutations each of which moves only finitely many elements. Then every $g\in G$ centralizes all but finitely many elements of each conjugacy class $h^G$.  

Given $g, h \in G\setminus \{e\}$, let us choose $C>0$ such that the finite subset of $X$ consisting of the elements moved by $g $  and the elements moved by $h$ has all elements within distance 
$\leq C $ each other. Since elements of $G$ are distance-preserving, for every conjugate $h^f$, $f\in G$, the elements moved by $h^f$ are also within distance $\leq C$ of each other. Hence, if any of the elements moved by $h^f$ has distance $>2C$ from each element moved by $g$, then the set of  elements moved by $h^f$  must be disjoint from the set moved by $g$, so  $h^f$ and $g$ commute. 
So, if $h^f$ and $g$ do not commute, the elements  moved by $h^f$ must lie within distance $\leq 2C$ of an arbitrary chosen element $x$  moved by $g$.
But the number of elements lying within that distance of $x$ if finite, so there are only finitely many posibilities for the permutation  
$h^f$.
\vskip 10pt

{\it Claim 2.} For $X$ and $G$ as in   Claim 1, if $X$ is infinite and $G$ is transitive on $X$, then every nonidentify element $g\in G$ has infinite conjugacy class $g^G$. 

Given finitely many conjugates $g_1, \dots , g_n$ of
$g$, we shall find another. 
Let $Y$ by the finite subset of $X$ consisting of all elements moved by any $g_1, \dots , g_n$, and again choose $C>0$ such that the distances between the element of $Y$ are all $\leq C$.
Since $X$ is infinite, the hypothesis of 
Claim 1 imply that distances among points of $X$ are unbounded, so as $G$ is transitive on $X$, we can find $h\in G$ carries a point moved by $g$ to a point at distance $>2C$ from point of $Y$. Hence, the set of point moved by $g^h$, namely, the translate by $h$ of the set moved by $g$, is not contained in $Y$, so $g^h \notin \{g_1, \dots , g_n \}$.
So, the conjugacy class of $g$ is indeed infinite.

\vskip 10pt

It remains to give an example of $X$ and $G$ with above properties. 

Let $X$ be the set of all  sequences $(a_1, a_2, \dots )$ of $0$'s and $1$'s such that almost all the $a_i$ are $0$.
Metrize $X$ by letting $d((a_1, a_2, \dots), (b_1, b_2, \dots ))$ be the greatest $n$ such  that $a_n \neq b_n,$ or $0$ if  
$(a_1, a_2, \dots)= (b_1, b_2, \dots )$.
That there are only finitely many elements distances $C$ of any element of $X$ is clear.

Let $G$ be the group of all distance-preserving permutations of $X$ which move only finitely many elements. 
We shall show that $G$ is transitive by constructing, for any $(a_1, a_2, \dots)\in X$
an element $g\in G$ which carries $(0,0, \dots )$ to $(a_1, a_2, \dots)$. Choose $n$ such that $a_i =0$ for all $i>n$. Let $g$ carries each element $(b_1, b_2, \dots )$ which likewise has $b_i =0$ for all $i>n$ to $(b_1 + a_1, b_2 +a_2, \dots )$, while fixing all other elements $(b_1, b_2, \dots )$. 
The verification of $g\in G$, and that $g$ carries 
$(0,0,\dots)$  to $(a_1, a_2, \dots)$ are straightforward. $ \  \Box$

\vskip 10pt

G. Bergman noticed that the group $G$ constructed in the proof of Theorem  4 can be  described as the direct limit $G_0 \longrightarrow G_1 \longrightarrow \dots \longrightarrow G_n \longrightarrow\dots$, where $G_0$ is trivial and $G_{n+1}= (G_n \times G_n) \leftthreetimes \mathbb{Z}_2$, with $\mathbb{Z}_2$
acting on $G_n \times G_n$ by interchanging the two coordinates, and with  $G_n$
embedded in $G_{n+1}$ by sending $g$ to $((g,e), e)$.

\vskip 7pt

We show that the answer to our question  is affirmative provided that $G$ is finitely generated. 
Let $F$ be a finite subset of $G$ such that $F= F^{-1}$, $e\in F$, $e$ is the identity of $G$ and $F$ generates $G$. 
We assume that each connected component of $\stackrel{\leftrightarrow}{G} $ is discrete, take an arbitrary element $g\in G$ and show that $g^G$ is finite. We act on $g$ by conjugations from $x\in F$, write each $g^x$ as a word in $F$ of minimal length, delete duplicates (i.e. words which define the same elements) and get a subset $A_0$. Then we repeat this procedure for each element $g\in A_0$ and get a subset $A_1$,  $A_0 \subseteq A_1$. Since $F$ is finite, by the assumption there exists $n\in \mathbb{N}$ such that $A_{n+1}= A_n$.
This means that $g^G = A_n$.

\section{Cellularity}

A coarse space $(X,\mathcal{E})$ is called {\it cellular} if  $\mathcal{E}$ has a base consisting of equivalence relations. By [12, Theorem 3.1.3], $(X,\mathcal{E})$ is cellular if and only if $asdim \ (X,\mathcal{E})=0$.

Applying Theorem 3.1.2 from \cite{b12}
we get 

\vskip 7pt

 (1)  $  \  \stackrel{\leftrightarrow}{G} $ is cellular if and only if, for every finitely generated subgroup $H$ of $G$, 
 
 there exists a finite subset $F$ of $G$ such that
$g^H \subseteq g^F$ for each $g\in G$.

\vskip 10pt

We recall that a group $G$ is {\it locally normal} if each finite subset of $G$ is contained in some finite normal subgroup and use the following characterization \cite{b2}

\vskip 10pt

(2)  $ \ G$ is an FC-group if and only if $G / Z_G$ is locally normal and each element 

of $G$ is contained in finitely generated normal subgroup, $Z_G$ is the center of $G$.

\vskip 10pt

A group $G$ is called {\it locally finite} if each finite subset of $G$ generates a finite subgroup. 

\vskip 10pt

{\bf Theorem 5.} {\it For a group $G$, 
$\stackrel{\leftrightarrow}{G} $ 
is cellular if and only if 
$G / Z_G$ is locally finite.

\vskip 10pt

Proof.} We suppose that $\stackrel{\leftrightarrow}{G} $ is 
cellular and show 

\vskip 7pt

(3) for every element $a\in G$ of infinite order there exists $n\in \mathbb{N}$ such that 

$a^n \in  Z_G$.

\vskip 10pt

We denote by $A$ the subgroup of $G$ generated by $a$ and use (1) to choose a finite subset $F$ of $G$ such that $g^A \subseteq g^F$ for each $g\in G$. Let $|F|=n$.
Since  $|g^A|\leq n $, $a^k g =ga^k$ for some $k\leq m$.
We put $n=m !$.

By (1), every finitely generated subgroup $H$ of $G$ is an FC-group. By (3), $H / (H\cap Z_G)$ is a torsion group. 
Applying (2), we conclude that $H / (H\cap Z_G)$ is finite. Hence, $G / Z_G$ is locally finite.

Now let $G / Z_G$ is locally finite. We take an arbitrary finitely generated subgroup $H$ of $G$, choose a set $h_1, \dots, h_n$ of representatives of right cosets of $H$ by 
$H\cap Z_G$, put $F=\{ h_1, \dots ,h_n \}$
and note that $g^H= g^F$ for each $g\in G$. 
Applying (1), we conclude that  $\stackrel{\leftrightarrow}{G} $ 
is cellular. $ \ \Box$

\vskip 10pt

{\bf Remark 6.}  
Every finitely generated subgroup of a group $G$ is an FC-group if and only if $g^H$  is finite for each $g\in G$ and every finitely generated subgroup $H$.
If $G / Z_G$ is locally finite then every finitely generated subgroup $H$ of $G$ is an FC-group. 
We show that the converse statement does not hold. Let 
$H= \oplus _{i<\omega} H_i$ be the direct sum of $\omega$
copies of $\mathbb{Z}_2$. We partition 
$\omega$ into consecutive intervals $\{W_i : i<\omega  \}$ of length $|W_i|=i+1$.  Then we take an automorphism $a$ of $H$ acting on each 
$\oplus \{H_m : m\in W_i \}$ as the cyclic permutations  of  coordinates, denote by $A$ the cyclic group generated by $A$ and consider the semidirect product $G=H \leftthreetimes A$.
Then every finitely generated subgroup of $G$ is an FC-group 
but $a^n \notin Z_G$ for each $n\in \mathbb{N}$ so 
$G/ Z_G$ is not locally finite.

\section{Asymorphic embeddings}

Let $(X, \mathcal{E})$, $(X^\prime , \mathcal{E}^\prime )$ be coarse spaces. A mapping $f: X\longrightarrow X^\prime$
is called {\it macro-uniform} if, for every  
$E\in \mathcal{E}$,  there exists $E^\prime \in \mathcal{E}^\prime $ such that 
$f(E[x]) \subseteq E^\prime [f(x)]$ for each $x\in X$.
We say that an injective mapping $f: X\longrightarrow X^\prime $
is an {\it asymorphic embedding} if 
$f: X\longrightarrow X^\prime $ and 
$f^{-1}: f(X)\longrightarrow X $ are macro-uniform.

\vskip 10pt

{\bf Theorem 7.} {\it Every finitary coarse space $(X, \mathcal{E})$ admits an asymorphic embedding to $\stackrel{\leftrightarrow}{G} $  for an appropriate choice of a group $G$.

\vskip 10pt

Proof.} We represent $(X, \mathcal{E})$ as the coarse space $X_H$  for some group $H$
 of permutations of $X$, see [7, Theorem 1]. 
 We consider $\{0, 1  \}^X$ as a  group with point-wise 
 addition $mod \ 2$.
 For $h\in H$  and $\chi\in \{ 0,1\}^X$, we put  $\chi_h (y)=\chi (h^{-1} y)$. 
 Then we define a semidirect product $G= \{ 0,1\}^X \leftthreetimes   H$ by 
 $$(\chi, h) (\chi^\prime, h ^\prime)=(\chi + \chi^\prime _h , hh^\prime) $$
 and note that the mapping $f: X\longrightarrow 
 \{ 0,1\}^X$, $f(x)$ is the characteristic function of $\{x\}$ is an asymorphic embedding of $(X, \mathcal{E})$ into 
 $\stackrel{\leftrightarrow}{G} $. $ \ \Box$
 
 \vskip 10pt
 
 If a subset $A$ of a coarse space $(X, \mathcal{E})$ is the union of $n$ discrete subsets then $A$ is $n$-discrete. 
 
 \vskip 10pt

{\bf Theorem 8.} {\it Let $G$ be a countable group. Then every $n$-discrete subset $A$ of 
$\stackrel{\leftrightarrow}{G} $  can be partitioned into $n $ discrete subsets. 

\vskip 10pt

Proof.} Use arguments proving this statement in the case of a connected coarse space with a linearly ordered base [6, Theorem 1.2]. 
$ \ \Box$

\vskip 10pt

{\bf Theorem 9.} {\it There exists a group $G$ such that 
$\stackrel{\leftrightarrow}{G} $  has 2-discrete subset which cannot be finitely partitioned into discrete subsets. 

\vskip 10pt

Proof.} By Theorem 6.3 from \cite{b3}, there exists 2-discrete finitary coarse  space on $\omega$ which  cannot be finitely partitioned into discrete subspaces. Apply Theorem 7.
$ \ \Box$

\section{The space of subgroups}

 \vskip 10pt
 
 For a group $G$ we denote by $\mathcal{S}$
 $(\stackrel{\leftrightarrow}{G} )$ the set 
$\mathcal{S} (G)$ of all subgroups of $G$ endowed with the coarse structure with the base 
$$\{\{ (X, Y) \in \mathcal{S} (G)\times  \mathcal{S} (G) : Y\in X^F \} : F\in [G]^{<\omega}\},$$
$X^F = \{ g^{-1} X g: g\in F\}.$

\vskip 10pt

We recall that $G$ is a {\it Dedekind group} if each subgroup of $G$ is normal. A non-abelian Dedekind group is called Hamiltonian. By \cite{b1}, 
\vskip 7pt

$ \ \ \ (4)$ $G$ is Hamiltonian if and only if $G$  is isomorphic to $Q_8 \times P$, where $Q_8$ is 

$ \ \ \ \ $ the quaternion group, $P$
 is an Abelian group without of elements of order 4.
 
 \vskip 10pt

{\bf Theorem 10.} {\it For an infinite group $G$, $\mathcal{S} (\stackrel{\leftrightarrow}{G} )$ is discrete if and only if $G$ is a Dedekind group.
 
\vskip 10pt

Proof.} If each subgroup of $G$ is normal then, evidently, $\mathcal{S} (G)$ is discrete.

We assume that $\mathcal{S} (\stackrel{\leftrightarrow}{G} )$ is discrete and consider two cases. 

\vskip 10pt

{\it Case 1:} $G$ has an element of infinite order. 
First, we show that every infinite cyclic subgroup of $G$ is invariant. We suppose the contrary and choose an infinite cyclic  subgroup $A, \ A=<a>$ and $z\in G$ such that $z^{-1} az \notin A$. Since 
$\mathcal{S}$
 $(\stackrel{\leftrightarrow}{G} )$ is discrete, there exists $m\in \mathbb{N}$ such that $z^{-1}<a^n>z = <a^n>$ for  each $n>m$.
 By the same reason, there exists $k\in \mathbb{N}$ such that $z^{-1} <aa^n>z = <aa^n>$ for each $n>k$.
 We take an arbitrary $n$ such that $n>m$, $n>k$. 
 Then
 $z^{-1} a^{n+1}
 z =(z^{-1} a z) (z^{-1} a^n z) \in <a^{n+1}> $, $z^{-1} a^{n} z \in <a^n>$, so $z^{-1} a^{n} z \in A$, contradicting the choice of $A$ and $z$. 
 
 Second, we take an arbitrary element $a\in G$ of infinite order and show that $a\in Z_G$. 
 Assuming the contrary, we  get $z\in G$ such that  $z^{-1} a z \neq a$.
 By above paragraph  $z^{-1} a z = a^{-1}$, so $z^{-2} a z^2 = a$ and 
 $(a^n z) (a^n z)= a^n z^2   z^{-1} a^n z =
  a^n z^2  a ^{-n} = z^2$
for each $n\in \mathbb{N}$. Since  
  $\mathcal{S} (\stackrel{\leftrightarrow}{G} )$ is discrete, there exists $m\in \mathbb{N}$ such that 
$$ z^{-1} (<a^n z > <z^2>) z = < a^n z > <z^2>$$  for each $n>m$. Hence, 
 $$ z^{-1} (a^n z ) z = a^{-n}z \in  < a^n z > <z^2>$$  and $a^{2n}\in <z>$, contradicting $z^{-1} a^{2n} z = a^{-2n} $.

 If $b$ is an element of finite order and $a$ is an element of infinite order then $ab$ has
 an infinite order because $a\in Z_G$, so $ab\in Z_G$, $b\in Z_G$, and $G$ is Abelian.   
\vskip 10pt

{\it Case 2:} Every element of $G$ has a finite order. We prove that $G$ is a Dedekind group provided that the following condition holds

\vskip 10pt
$ \ \ (5) $  for every finite subset $K$ of $G$ containing the identity $e$, there exists $a\in G$, 

$ \ \ \ \ $ $a\neq e$ such that $K\cap <a>=\{e\}$.

\vskip 10pt

We suppose the contrary and choose $b\in G$, $z\in G$ such that $z^{-1}bz \notin <b>$. Since 
$\mathcal{S} (\stackrel{\leftrightarrow}{G} )$ is discrete, by (5), there exists $a\in G$, $a\neq e$ such that 
$$z^{-1}bz  <b> \cap <a> =\{e\}, \ z^{-1} <a>z = <a>, $$
$$b^{-1}<a> b = <a>,  \  z^{-1}<b><a>z = <b><a>.$$
Then 
$z^{-1}baz =(z^{-1}bz) (z^{-1}az) \in <b><a>$, 
$z^{-1}bz \in <b><a>$ and $z^{-1}bz \in <b>$, contradicting the choice of $b$  and $z$.

We denote by $\pi (G)$ the set of all prime divisors 
 of orders of elements of $G$ and put $X_n = \{ g\in G: g^n = e\}$.
 If $G$ is not a Dedekind group, by (5), $\pi (G)$  is finite and $X_p$ is finite for each $p\in \pi (G)$. We prove that $G$ is layer-finite: $X_n$ is finite for each $n\in \mathbb{N}$. 
 It suffices to verify that $X_{p^n} $  is finite for all $p\in \pi (G)$, $n\in \mathbb{N}$. We suppose that $X_{p^m} $ is finite but $X_{p^{m+1}} $
 is infinite.
 Then there exists a sequence $(a_n)_{n\in\omega}$ in $G$ and $a\in G$ such that $|a_n|= p^{m+1}$, $|a|=p^m$ and 
 $<a_n>\cap <a_k>=<a>$ for all distinct $n,k\in \mathbb{N}$. We denote by $H$ the subgroup of $G$ generated by the set  $\{ a_n: n\in \omega \} $ and put $M= H/<a>$. 
 Since $\mathcal{S}$
$(\stackrel{\leftrightarrow}{M} )$ is discrete, applying (5) and (4) to $M$, we conclude that $M$ has an infinite Abelian subgroup of exponent $p$.
 By the Gr$\stackrel{..}{u}$n's lemma (see [5], p. 398),  $H$ has an infinite Abelian subgroup of exponent $p$,  so $X_p$ is infinite and we get a contradiction.

Thus, our assumption that $G$ is not a Dedekind group gives $G$ is layer-finite and $\pi (G)$ is finite.
Since $G$ is infinite, by the Chernikov's theorem [4], $G$ has a central quasi-cyclic $p$-group $A$, 
$A=\cup _{n\in\omega} <a_n>, \ a_{n+1}^p =a_n.$
We take $c, z\in G$  such that $z^{-1} c z \neq <c>$, $|c|=q^m$,  $q\in \pi (G)$.
Since $\mathcal{S}$ $(\stackrel{\leftrightarrow}{G} )$ is discrete, there exists $k\in \mathbb{N}$ such that, for each $n>k$, we have  
$$z^{-1}<a_n c> z = <a_n c>,  \ \  a_n (z^{-1}c z)\in <a_n c>.$$
If $q\neq p$ then $z^{-1} c z\in <c>$, contradicting the choice of $c$ and $z$.
If $q=p$ and $n> 2m$, $n>k$ then 
$(a_n c)^{p^m}= a_n ^{p^m} $, $|a_n ^{p^m}|> p^m $
and $z^{-1} c z \in <a_n ^{p^m}>$.
Since $A$
is central, $z^{-1} c z =c$ and $z^{-1} c z \in <c>$, contradicting the choice of $z, c$.
The proof is completed. $ \ \Box$

\vskip 10pt

{\bf Remark 11.} Let $G$ be a transitive group of permutations of a set $X$, $St (x) = \{ g\in G: gx = x\}$, $x\in X$. 
Then the natural mapping $x\mapsto  St(x)$ is an asymorphic embedding of the finitary coarse space $X_G$ into $\mathcal{S}$
 $(\stackrel{\leftrightarrow}{G} )$.

 \vskip 10pt
 
 If 
$(\stackrel{\leftrightarrow}{G} )$ is cellular then applying (1) we see that 
 $\mathcal{S}$
 $(\stackrel{\leftrightarrow}{G} )$ is cellular. 

\vskip 10pt

{\bf Question 12.} Is $\stackrel{\leftrightarrow}{G}$ cellular provided that 
 $\mathcal{S}$
 $(\stackrel{\leftrightarrow}{G} )$ is cellular? 

\vskip 15pt

\section{The direct union of connected components }

Let $(X, \mathcal{E} )$ be a coarse space, 
$\lbrace X_\alpha : \alpha < \kappa\rbrace$ is the set of all connected components of 
$(X, \mathcal{E} )$. We say that $(X, \mathcal{E} )$ is the {\it direct union} of $\lbrace X_\alpha :\alpha<\kappa \rbrace $ if, for each 
$E\in \mathcal{E}$, there exists  $\alpha_1, \dots, \alpha_n$ such that $E[x]=\lbrace x\rbrace$ for each 
$x\in X_\alpha$, 
$\alpha< \kappa$, $\alpha\notin \lbrace \alpha_1, \dots, \alpha_n\rbrace$.

If a group $G$ is either Abelian or the set of conjugacy classes of $G$ is finite then 
$\stackrel{\leftrightarrow}{G}$ is the direct union of conjugacy classes.

For every natural number $n$, G. Bergman used $HNN$-extensions to construct a group $G$ such that $G$ has an infinite center (so the number of conjugacy classes of $G$ is infinite)  and only $n$ conjugacy classes of $G$ are not singletons. 
Also, he proved that if $\stackrel{\leftrightarrow}{G}$ is the direct union of  
conjugacy classes then all but finely many  conjugacy classes  are singletons. 
\vskip 10pt

{\bf Acknowledgment}. We thank George Bergman for the kind suggestion to incorporate his results in this paper.


\vskip 15pt

CONTACT INFORMATION
\vskip 15pt

I.~Protasov: \\
Faculty of Computer Science and Cybernetics  \\
        Taras Shevchenko National University of Kyiv \\
         Academic Glushkov pr. 4d  \\
         03680 Kyiv, Ukraine \\ i.v.protasov@gmail.com

K.~Protasova: \\
Faculty of Computer Science and Cybernetics  \\
Taras Shevchenko National University of Kyiv \\
         Academic Glushkov pr. 4d  \\
         03680 Kyiv, Ukraine \\ k.d.ushakova@gmail.com

\end{document}